\documentclass[12pt,a4paper,oneside,reqno,notitlepage]{amsart}
\usepackage{amssymb}
\textwidth
150mm

\theoremstyle{plain}

\numberwithin{equation}{section}

\begin{document}

\baselineskip 8mm
\parindent 9mm

\title[]
{Square-mean weighted pseudo almost automorphic solutions for stochastic semilinear integral equations}

\author{Kexue Li, \ \   Jigen Peng }

\address{School of Mathematics and Statistics, Xi'an
Jiaotong University, Xi'an 710049, China}
\email{ kxli@mail.xjtu.edu.cn,\ jgpeng@mail.xjtu.edu.cn}

\thanks{{\it 2010 Mathematics Subjects Classification}: 34F05}

\thanks{This work was supported by the Natural Science Foundation of
China under the contact No. 11201366 and 11131006}

\begin{abstract}
In this paper, we introduce the concept of $S^{2}$-weighted pseudo almost automorphy for stochastic processes. We study the existence and uniqueness of square-mean weighted pseudo almost automorphic solutions for the semilinear stochastic integral equation
$x(t)=\int_{-\infty}^{t}a(t-s)[Ax(s)+f(s,x(s))]ds+\int_{-\infty}^{t}a(t-s)\varphi(s,x(s))dw(s), \ t\in\mathbb{R}$, where $a\in L^{1}(\mathbb{R}_{+})$, $A$ is the generator of an integral resolvent family on a Hilbert space $H$, $w(t)$ is the two-sided $Q$-Wiener process, $f,\varphi: \mathbb{R}\times L^{2}(P,H)\rightarrow L^{2}(P,H)$ are two $S^{2}$-weighted pseudo almost automorphic functions.
\end{abstract}
\maketitle

\section{\textbf{Introduction}}

The almost automorphic function introduced by Bochner \cite{Bochner} is a generation of the almost periodic function. Henriquez and Lizama \cite{HC} investigated the existence and regularity of compact almost automorphic solutions for the semilinear integral equation
\begin{align}\label{integral}
u(t)=\int_{-\infty}^{t}a(t-s)[Au(s)+f(s,u(s))]ds, \ t\in \mathbb{R},
\end{align}
where $A: D(A)\subset X\rightarrow X$ is the generator of an integral resolvent defined on a Banach space $X$, $a\in L^{1}(\mathbb{R}_{+})$.

A rich source for vector-valued Volterra equations is the continuum in mechanics for materials with memory, i.e. the theory of viscoelastic materials (see \cite{Pruss}). Cuevas and Lizama \cite{CC} studied the existence and uniqueness of almost automorphic solutions for problem (\ref{integral}). Zhang \cite{Zhang} introduced the pseudo almost periodic function and studied the pseudo almost periodic solutions of some differential equations. Liang \cite{JJT} developed the theory of the pseudo almost automorphic function, which is a generation of the pseudo almost periodic function. Zhao et al \cite{ZYG} considered the existence and uniqueness of pseudo-almost automorphic mild solutions for problem (\ref{integral}).

Diagana and   N'Gu$\acute{e}$r$\acute{e}$kata \cite{DG} introduced the concept of  $S^{p}$-almost automorphy to study the existence of solutions to some semilinear equations. The $S^{p}$-almost automorphic function is a generalization of the almost automorphic function. Ding et al. \cite{Ding} established a composition theorem about $S^{p}$-almost automorphic functions and studied the existence and uniqueness of almost automorphic solutions for semilinear evolution equations in Banach space. Lizama and Ponce \cite{CR} studied the existence of an almost automorphic mild solution for problem (\ref{integral}) with $S^{p}$-almost automorphic coefficients.  Zhang et al. \cite{RYG} introduced the concept of $S^{p}$-weighted pseudo almost automorphy and studied the existence and uniqueness of  $S^{p}$-weighted pseudo almost automorphic solutions to nonautomous evolution equations.  Xia and Fan \cite{XF} proposed the concept of weighted Stepanov-like pseudo
almost automorphic functions,  studied the properties of the
weighted Stepanov-like pseudo almost automorphic functions in Banach space. Zhang and Chang \cite{RY} investigated the existence of weighted pseudo almost automorphic solution for problem (\ref{integral}) with $S^{p}$-weighted pseudo almost automorphic coefficients.

In the real world, stochastic noises are unavoidable. Fu and Liu  \cite{FL} introduced the square-mean almost automorphic process and investigated the square-mean almost automorphic mild solution for a class of stochastic differential equations. Chen and Lin  \cite{CL} introduced the concept of the square-mean pseudo almost automorphic for a stochastic process and studied the existence, uniqueness and global stability of the square-mean pseudo almost automorphic solutions for a class of stochastic evolution equations. Chen and Lin \cite{Chenzhang} introduced the concepts and properties of the square-mean weighted pseudo almost automorphy for a stochastic process, they investigated the well-posedness of the square-mean weighted pseudo almost automorphic solutions for a class of non-autonomous stochastic differential equations. Chang et al. introduced the concept of Stepanov-like almost automorphy (or $S^{2}$-almost automorphy) for stochastic processes. They used the results obtained to investigate
the existence and uniqueness of a Stepanov-like almost automorphic mild solution to a
class of nonlinear stochastic differential equations in a real separable Hilbert space.

To the best of our knowledge, there are no considerations for problem (\ref{integral}) perturbed by stochastic noises. In this paper, we study the existence of weighted pseudo almost automorphic solutions for following the stochastic integral equation
\begin{align}\label{key}
x(t)=\int_{-\infty}^{t}a(t-s)[Ax(s)+f(s,x(s))]ds+\int_{-\infty}^{t}a(t-s)\varphi(s,x(s))dw(s), \ t\in\mathbb{R},
\end{align}
where $a\in L^{1}(\mathbb{R}_{+})$, $A$ is the generator of an integral resolvent family on a Hilbert space $H$, $w(t)$ is the two-sided $Q$-Wiener process, $f,\varphi: \mathbb{R}\times L^{2}(P,H)\rightarrow L^{2}(P,H)$ are two $S^{2}$-weighted pseudo almost automorphic functions satisfying some conditions.

The paper is organized as follows. In Section 2, we recall some concepts and preliminary facts. In Section 3, we  introduce the concept of  $S^{2}$-weighted pseudo almost automorphy and obtain some properties of $S^{2}$-weighted pseudo almost automorphic functions. In Section 4, we study the existence and uniqueness of square-mean weighted pseudo almost automorphic solutions for problem (\ref{key}).

\section{Preliminaries}
Let $(H,\|\cdot\|)$ be a separable real Hilbert space and let $(\Omega, \mathcal{F}, P)$ be a complete probability space. $L^{2}(P,H)$ denotes the space of all $H$-valued random variables $x$ such that $\mathbb{E}\|x\|^{2}=\int_{\Omega}\|x\|^{2}d P<\infty$, which is a Banach space equipped with the norm
$\|x\|_{2}=(\mathbb{E}\|x\|^{2})^{\frac{1}{2}}$. For $f: \mathbb{R}\rightarrow \mathbb{R}$, we denote $\lim_{t\rightarrow \infty}\sup f(t):=\inf_{M>0}\big(\sup_{|t|>M}f(t)\big)$.  \\
\textbf{Definition 2.1.} A stochastic process $x: \mathbb{R}\rightarrow L^{2}(P,H)$ is said to be stochastically continuous if
\begin{align*}
\lim_{t\rightarrow s}\mathbb{E}\|x(t)-x(s)\|^{2}=0.
\end{align*}
\textbf{Definition 2.2.} A stochastic process $x: \mathbb{R}\rightarrow L^{2}(P,H)$ is said to be stochastically bounded if there exists a constant $M>0$ such that
\begin{align*}
\mathbb{E}\|x(t)\|^{2}\leq M.
\end{align*}

Denote by $SBC(\mathbb{R}, L^{2}(P,H))$ the collection of all stochastically continuous  and bounded processes.

Let $\mathcal{U}$ be the set of all functions that are positive and locally integrable over $\mathbb{R}$. For given $T>0$ and $\rho\in \mathcal{U}$, define $\mu(T,\rho)=\int_{-T}^{T}\rho(t)dt$, $\mathcal{U}_{\infty}=\{\rho\in \mathcal{U}: \lim_{T\rightarrow \infty}\mu(T,\rho)=+\infty\}$, and
$\mathcal{U}_{b}=\{\rho\in \mathcal{U}_{\infty}: \rho \ \mbox{is bounded and}\ \inf_{t\in\mathbb{R}}\rho(t)>0\}$. \\
\textbf{Definition 2.3.} For $\rho_{1}, \rho_{2}\in \mathcal{U}_{\infty}$, if $\frac{\rho_{1}}{\rho_{2}}\in \mathcal{U}_{b}$, then we call $\rho_{1}$ is equivalent to $\rho_{2}$, and write $\rho_{1}\sim\rho_{2}$.

For  $\rho\in \mathcal{U}_{\infty}$, $s\in \mathbb{R}$, define $\rho_{s}(t)=\rho(t+s), t\in \mathbb{R}$  and  $\mathcal{U}_{T}=\{\rho\in\mathcal{U}_{\infty}: \rho\sim \rho_{s}, \forall s\in\mathbb{R}\}$.\\
\textbf{Definition 2.4.} For $\rho\in \mathcal{U}_{\infty}$, define
\begin{align*}
SPAA_{0}(\mathbb{R}, L^{2}(P,H),\rho)&=\{x\in SBC(\mathbb{R}, L^{2}(P,H)): \lim_{T\rightarrow\infty}\frac{1}{\mu(T,\rho)}\int_{-T}^{T}\left(\mathbb{E}\|x(t)\|^{2}\right)^{\frac{1}{2}}\rho(t)dt=0\}
\end{align*}
and
\begin{align*}
&SPAA_{0}(\mathbb{R}\times L^{2}(P,H),L^{2}(P,H),\rho)=\{f:\mathbb{R}\times L^{2}(P,H)\rightarrow L^{2}(P,H) \ \mbox{is stochastically continuous}, \\
\ &f(\cdot,x) \ \mbox{is stochastically bounded for each}\  x\in L^{2}(P,H) \ \mbox{and} \\ &\lim_{T\rightarrow\infty}\frac{1}{\mu(T,\rho)}\int_{-T}^{T}\left(\mathbb{E}\|f(t,x)\|^{2}\right)^{\frac{1}{2}}\rho(t)dt=0\}.
\end{align*}
\textbf{Definition 2.5.} A subset $\mathcal{D}$ of $SBC(\mathbb{R}, L^{2}(P,H))$ is said to be translation invariant if for any $x(t)\in \mathcal{D}$, $x(t+\tau)\in \mathcal{D}$ for any $\tau\in\mathbb{R}$.

For simplicity, denote $\mathcal{U}^{inv}=\{\rho\in \mathcal{U}_{\infty}: SPAA_{0}(\mathbb{R}, L^{2}(P,H),\rho) \ \mbox{is translation invariant}\}$.
\textbf{Definition 2.6.} (\cite{FL}) A stochastically continuous stochastic process $x: \mathbb{R}\rightarrow L^{2}(P,H)$ is said to be square-mean almost automorphic if every sequence of real numbers $\{s'_{n}\}$ has a subsequence $\{s_{n}\}$ such that for some stochastic process $y: \mathbb{R}\rightarrow L^{2}(P,H)$
\begin{align*}
\lim_{n\rightarrow \infty}\mathbb{E}\|x(t+s_{n})-y(t)\|^{2}=0 \ \mbox{and}\ \lim_{n\rightarrow \infty}\mathbb{E}\|y(t-s_{n})-x(t)\|^{2}=0
\end{align*}
holds for each $t\in \mathbb{R}$. The collection of all square-mean almost automorphic stochastic processes $x: \mathbb{R}\rightarrow L^{2}(P,H)$ is denoted by $SAA(\mathbb{R}, L^{2}(P,H))$. \\
\textbf{Lemma 2.7.} (\cite{FL}) $SAA(\mathbb{R}, L^{2}(P,H))$ is a Banach space when it is equipped with the norm $\|x\|_{\infty}= \sup_{t\in \mathbb{R}}(\mathbb{E}\|x(t)\|^{2})^{\frac{1}{2}}$. \\
\textbf{Definition 2.8.}  Let $\rho\in \mathbb{U}_{\infty}$.  A function $f\in SBC(\mathbb{R}, L^{2}(P,H))$ is called square-mean weighted pseudo
almost automorphic if it can be expressed as $f=g+h$, where $g\in SAA(\mathbb{R}, L^{2}(P,H))$ and $h\in SPAA_{0}(\mathbb{R}, L^{2}(P,H),\rho)$. The collection of all such processes will be denoted by $WPAA(\mathbb{R}, L^{2}(P,H),\rho)$. \\
\textbf{Lemma 2.9.} Suppose $\rho\in \mathcal{U}^{inv}$ and $f=g+\varphi\in WPAA(\mathbb{R}, L^{2}(P,H),\rho)$, where $g\in SAA(\mathbb{R}, L^{2}(P,H))$, $\varphi(t)\in SPAA_{0}(\mathbb{R}, L^{2}(P,H),\rho)$. Then $\overline{\{f(t)|t\in \mathbb{R}\}}\supset \{g(t)|t\in \mathbb{R}\}$.\\
\textbf{Proof.} We show this lemma by contradiction. Assume that there exist $t_{0}\in R$  and $\varepsilon>0$ such that
\begin{align*}
\left(\mathbb{E}\|g(t_{0})-f(t)\|^{2}\right)^{\frac{1}{2}}\geq 2\varepsilon, \ t\in \mathbb{R}.
\end{align*}
Denote
\begin{align*}
B_{\varepsilon}=\{\tau\in\mathbb{R}|\left(\mathbb{E}\|g(t_{0}+\tau)-g(t_{0})\|^{2}\right)^{\frac{1}{2}}<\varepsilon\}.
\end{align*}
By Lemma 2.1 in \cite{CL}, there exist $s_{1},\ldots, s_{m}\in \mathbb{R}$ such that $\bigcup_{i=1}^{m}(s_{i}+B_{\varepsilon})=\mathbb{R}$.
For any $t\in \mathbb{R}$, there exists some $i$ such that $t\in s_{i}+B_{\varepsilon}$. By the Minkowski inequality,
\begin{align*}
(\mathbb{E}\|\varphi(t+s_{0}-s_{i})\|^{2})^{\frac{1}{2}}&\geq \left((\mathbb{E}\|f(t+t_{0}-s_{i})-g(t_{0})\|^{2})^{\frac{1}{2}}-(\mathbb{E}\|g(t_{0})-g(t+t_{0}-s_{i})\|^{2})^{\frac{1}{2}}\right)\\
&\geq \varepsilon.
\end{align*}
Then for $T>0$,
\begin{align*}
\frac{1}{\mu(T,\rho)}\int_{-T}^{T}(\mathbb{E}\|\varphi(t+s_{0}-s_{i})\|^{2})^{\frac{1}{2}}\rho(t)dt\geq \frac{\varepsilon}{\mu(T,\rho)}\int_{-T}^{T}\rho(t)dt=\varepsilon.
\end{align*}
On the other hand, since $\rho\in \mathcal{U}^{inv}$ and $\varphi(t)\in SPAA_{0}(\mathbb{R}, L^{2}(P,H),\rho)$ then
\begin{align*}
\lim_{T\rightarrow\infty}\frac{1}{\mu(T,\rho)}\int_{-T}^{T}(\mathbb{E}\|\varphi(t+s_{0}-s_{i})\|^{2})^{\frac{1}{2}}\rho(t)dt=0.
\end{align*}
This implies a contradiction. Therefore $\overline{\{f(t)|t\in \mathbb{R}\}}\supset \{g(t)|t\in \mathbb{R}\}$. The proof is complete. \ \ $\Box$ \\
\textbf{Theorem 2.10.}  If $\rho\in \mathcal{U}^{inv}$, then $(WPAA(\mathbb{R},L^{2}(P,H),\rho), \|\cdot\|_{\infty})$ is a Banach space.\\
\textbf{Proof.} The proof is similar to the proof of Theorem 2.5 in \cite{Blot} by minor revisions. So it is omitted. \ \ $\Box$ \\
\textbf{Definition 2.11.}(\cite{PT}) The Bochner transform $x^{b}(t,s)$, $t\in \mathbb{R}$, $s\in [0,1]$, of a stochastic process $x:\mathbb{R}\rightarrow L^{2}(P,H)$ is defined by
\begin{align*}
x^{b}(t,s)=x(t+s).
\end{align*}
\textbf{Definition 2.12.} The Bochner transform $f^{b}(t,s,x)$, $t\in \mathbb{R}$, $s\in [0,1]$, $x\in L^{2}(P,H)$, of a function $x:\mathbb{R}\times L^{2}(P,H)\rightarrow L^{2}(P,H)$ is defined by
\begin{align*}
f^{b}(t,s,x)=f(t+s,x)
\end{align*}
for each $x\in L^{2}(P,H)$.\\
\textbf{Definition 2.13.} (\cite{PT}) The space $BS^{2}(L^{2}(P,H))$ of all Stepanov bounded stochastic process consists of all stochastic processes $x$ on $\mathbb{R}$ with values in $L^{2}(P,H)$ such that $x^{b}\in L^{\infty}(\mathbb{R},L^{2}((0,1), L^{2}(P,H)))$. This is a Banach space with the norm
\begin{align*}
\|x\|_{S^{2}}=\|x^{b}\|_{L^{\infty}(\mathbb{R}, L^{2})}&=\sup_{t\in \mathbb{R}}\left(\int_{0}^{1}\mathbb{E}\|x(t+s)\|^{2}ds\right)^{\frac{1}{2}}\\
&=\sup_{t\in \mathbb{R}}\left(\int_{t}^{t+1}\mathbb{E}\|x(\tau)\|^{2}d\tau\right)^{\frac{1}{2}}.
\end{align*}
\textbf{Definition 2.14.}  (\cite{Chang}) A stochastic process $x\in BS^{2}(L^{2}(P,H))$ is called Stepanov-like almost automorphic (or $S^{2}$-almost automorphic) if $x^{b}\in SAA(\mathbb{R}, L^{2}((0,1), L^{2}(P,H)))$.  In other words, a stochastic process $x\in L^{2}_{loc}(\mathbb{R}, L^{2}(P,H))$ is said to be Stepanov-like almost automorphic if its Bochner transform $x^{b}: \mathbb{R}\rightarrow L^{2}((0,1), L^{2}(P,H))$ is square-mean almost automorphic in the sense that for every sequence of real numbers $\{s'_{n}\}$, there exist a subsequence $\{s_{n}\}$ and a stochastic process $y\in L^{2}_{loc}(\mathbb{R}, L^{2}(P,H))$ such that
\begin{align*}
\int_{t}^{t+1}\mathbb{E}\|x(s+s_{n})-y(s)\|^{2}ds\rightarrow 0 \ \mbox{and} \ \int_{t}^{t+1}\mathbb{E}\|y(s-s_{n})-x(s)\|^{2}ds\rightarrow 0
\end{align*}
as $n\rightarrow \infty$ pointwise on $\mathbb{R}$. The collection of all such processes will be denoted by $AS^{2}(\mathbb{R}, L^{2}(P,H))$. \\
\textbf{Lemma 2.15.} (\cite{Chang}) $AS^{2}(\mathbb{R}, L^{2}(P,H))$ is a Banach space when it is equipped with the norm $\|\cdot\|_{S^{2}}$. \\
\textbf{Definition 2.16.} (\cite{Chang}) A function $f: \mathbb{R}\times L^{2}(P,H)$, $(t,x)\rightarrow f(t,x)$ with $f(\cdot,x)\in L^{2}_{loc}(\mathbb{R}, L^{2}(P,H))$ for each $x\in L^{2}(P,H)$, is said to be $S^{2}$-almost automorphic in $t\in \mathbb{R}$ uniformly in $x\in L^{2}(P,H)$. That means, for every sequence of real numbers $\{s'_{n}\}$, there exist a subsequence $\{s_{n}\}$ and a function $\widetilde{f}: \mathbb{R}\times L^{2}(P,H)\rightarrow L^{2}(P,H)$ with $\widetilde{f}(\cdot,x)\in L^{2}_{loc}(\mathbb{R}, L^{2}(P,H))$ such that
\begin{align*}
\int_{t}^{t+1}\mathbb{E}\|f(s+s_{n},x)-\widetilde{f}(s,x)\|^{2}ds\rightarrow 0 \ \mbox{and} \ \int_{t}^{t+1}\mathbb{E}\|\widetilde{f}(s-s_{n},x)-f(s,x)\|^{2}ds\rightarrow 0
\end{align*}
as $n\rightarrow \infty$ pointwise on $\mathbb{R}$ and for each $x\in L^{2}(P,H)$. The collection of such functions will be denoted by  $AS^{2}(\mathbb{R}\times L^{2}(P,H), L^{2}(P,H))$. \\

\section{$S^{2}$-weighted pseudo almost automorphic stochastic processes}
In this section, we will introduce the concept of  $S^{2}$-weighted pseudo almost automorphy and obtain some properties of $S^{2}$-weighted pseudo almost automorphic functions. \\
\textbf{Definition 3.1.} Let $\rho\in \mathcal{U}_{\infty}$. A stochastic process $f\in BS^{2}(L^{2}(P,H))$ is said to be Stepanov-like weighted pseudo almost automorphic (or $S^{2}$-weighted pseudo almost automorphic) if it can be expressed as $f=g+h$, where $g\in AS^{2}(\mathbb{R}, L^{2}(P,H))$ and $h^{b}\in SPAA_{0}(\mathbb{R}, L^{2}((0,1), L^{2}(P,H),\rho))$, i.e.,
\begin{align*}
\lim_{T\rightarrow \infty}\frac{1}{\mu(T,\rho)}\int_{-T}^{T}\left(\int_{t}^{t+1}\mathbb{E}\|h(s)\|^{2}ds\right)^{\frac{1}{2}}\rho(t)dt=0.
\end{align*}
The collection of such functions will be denoted by $WPAAS^{2}(\mathbb{R}, L^{2}(P,H),\rho)$.\\
\textbf{Lemma 3.2.} Suppose $\rho\in \mathcal{U}^{inv}$. If $f\in WPAA(\mathbb{R}, L^{2}(P,H),\rho)$ then $f\in WPAAS^{2}(\mathbb{R}, L^{2}(P,H),\rho)$. \\
\textbf{Proof.} Let $f=g+h$, where $g\in SAA(\mathbb{R}, L^{2}(P,H))$ and $h\in SPAA_{0}(\mathbb{R}, L^{2}(P,H),\rho)$. It is easy to see that $g\in SAA(\mathbb{R}, L^{2}(P,H))\subset AS^{2}(\mathbb{R}, L^{2}(P,H))$. Next we show that $h^{b}\in SPAA_{0}(\mathbb{R}, L^{2}((0,1), L^{2}(P,H),\rho))$.
For $T>0$,
\begin{align*}
\int_{-T}^{T}\left(\int_{0}^{1}\mathbb{E}\|h(t+s)\|^{2}ds\right)^{\frac{1}{2}}\rho(t)dt&\leq \int_{-T}^{T}\left(\int_{0}^{1}\sup_{s\in[0,1]}\mathbb{E}\|h(t+s)\|^{2}ds\right)^{\frac{1}{2}}\rho(t)dt\\
&\leq \int_{-T}^{T}\left(\sup_{s\in[0,1]}\mathbb{E}\|h(t+s)\|^{2}\right)^{\frac{1}{2}}\rho(t)dt.
\end{align*}
Let $s_{0}\in [0,1]$ such that $\sup_{s\in [0,1]}\mathbb{E}\|h(t+s)\|^{2}=\mathbb{E}\|h(t+s_{0})\|^{2}$. Then we have
\begin{align*}
\frac{1}{\mu(T,\rho)}\int_{-T}^{T}\left(\int_{0}^{1}\mathbb{E}\|h(t+s)\|^{2}ds\right)^{\frac{1}{2}}\rho(t)dt&\leq \frac{1}{\mu(T,\rho)}\int_{-T}^{T}\left(\sup_{s\in [0,1]}\mathbb{E}\|h(t+s)\|^{2}\right)^{\frac{1}{2}}\rho(t)dt\\
&\leq\frac{1}{\mu(T,\rho)}\int_{-T}^{T}\left(\mathbb{E}\|h(t+s_{0})\|^{2}\right)^{\frac{1}{2}}\rho(t)dt.
\end{align*}
Since $\rho\in \mathcal{U}^{inv}$, we have
\begin{align*}
\lim_{T\rightarrow \infty}\frac{1}{\mu(T,\rho)}\int_{-T}^{T}\left(\mathbb{E}\|h(t+s_{0})\|^{2}\right)^{\frac{1}{2}}\rho(t)dt=0.
\end{align*}
This implies that $h^{b}\in SPAA_{0}(\mathbb{R}, L^{2}((0,1), L^{2}(P,H),\rho))$. The proof is complete. \ $\Box$ \\
\textbf{Lemma 3.3.} Suppose $\rho\in \mathcal{U}_{T}$ and $\lim_{T\rightarrow \infty}\sup{\frac{\mu(T-\eta, \rho)}{\mu(T,\rho)}}>0$ for every $\eta\in \mathbb{R}_{+}$. If $f=g+h\in WPAAS^{2}(\mathbb{R}, L^{2}(P,H),\rho)$, where $g\in AS^{2}(\mathbb{R}, L^{2}(P,H))$ and $h^{b}\in SPAA_{0}(\mathbb{R}, L^{2}((0,1), L^{2}(P,H),\rho))$. Then $\{g(t+\cdot):t\in \mathbb{R}\}\subset \overline{\{f(t+\cdot):t\in \mathbb{R}\}}$. \\
\textbf{Proof.} We prove this lemma by contradiction. If this is not true, then there are $t_{0}\in \mathbb{R}$ and $\varepsilon>0$ such that
\begin{align*}
\|g(t_{0}+\cdot)-f(t+\cdot)\|_{2}\geq 2\varepsilon, \ t\in \mathbb{R}.
\end{align*}
Let $B_{\varepsilon}:=\{\tau\in \mathbb{R}: \|g(t_{0}+\tau+\cdot)-g(t_{0}+\cdot)\|_{2}<\varepsilon\}$. Since $g^{b}\in SAA(\mathbb{R}, L^{2}((0,1), L^{2}(P,H)))$, similar to Lemma 2.1 in \cite{CL}, there exist $s_{1}, \ldots, s_{m}\in \mathbb{R}$ such that $\bigcup_{i=1}^{m}(s_{i}+B_{\varepsilon})=\mathbb{R}$. Let $\hat{s_{i}}=s_{i}-t_{0} (1\leq i\leq m)$, $\eta=\max_{1\leq i\leq m}|\hat{s_{i}}|$. For
$T\in \mathbb{R}$ with $|T|>\eta$, set $B^{(i)}_{\varepsilon, T}=[-T+\eta-\hat{s_{i}}, T-\eta-\hat{s_{i}}]\cap (t_{0}+B_{\varepsilon})$, $1\leq i\leq m$. We have $\bigcup_{i=1}^{m}\big(\hat{s_{i}}+B^{(i)}_{\varepsilon, T}\big)=[-T+\eta, T-\eta]$. Since $B^{(i)}_{\varepsilon, T}\subset [-T,T]\cap (t_{0}+B_{\varepsilon}), i=1,\ldots,m$, we have
\begin{align}\label{a}
\mu(T-\eta,\rho)&=\int_{-T+\eta}^{T-\eta}\rho(t)dt\nonumber\\
&\leq \sum_{i=1}^{m}\int_{\hat{s_{i}}+B^{(i)}_{\varepsilon, T}}\rho(t)dt\nonumber\\
&\leq \sum_{i=1}^{m}\int_{B^{(i)}_{\varepsilon, T}}\rho(t+\hat{s_{i}})dt\nonumber\\
&\leq \max_{1\leq i\leq m}\{a_{i}\}\sum_{i=1}^{m}\int_{[-T,T]\cap (t_{0}+B_{\varepsilon})}\rho(t)dt\nonumber\\
&=m\cdot\max_{1\leq i\leq m}\{a_{i}\}\int_{[-T,T]\cap (t_{0}+B_{\varepsilon})}\rho(t)dt,
\end{align}
where $a_{i}=\lim_{t\rightarrow \infty}\sup\frac{\rho(t+\widehat{s_{i}})}{\rho(t)}\in (0, \infty)$. \\
By the Minkowski inequality, for every $t\in t_{0}+B_{\varepsilon}$, we have
\begin{align*}
\|h(t+\cdot)\|_{2}&=\|f(t+\cdot)-g(t+\cdot)\|_{2}\\
&\geq \|g(t_{0}+\cdot)-f(t+\cdot)\|_{2}-\|g(t+\cdot)-g(t_{0}+\cdot)\|_{2}\\
&\geq \varepsilon.
\end{align*}
Then
\begin{align}\label{b}
\frac{1}{\mu(T,\rho)}\int_{-T}^{T}\left(\int_{t}^{t+1}\mathbb{E}\|h(s)\|^{2}ds\right)^{\frac{1}{2}}\rho(t)dt&=\frac{1}{\mu(T,\rho)}\int_{-T}^{T}\rho(t)\|h(t+\cdot)\|_{2}dt\nonumber\\
& \geq \frac{1}{\mu(T,\rho)}\int_{[-T,T]\cap(t_{0}+B_{\varepsilon})}\rho(t)\|h(t+\cdot)\|_{2}dt\nonumber\\
& \geq \frac{\varepsilon}{\mu(T,\rho)}\int_{[-T,T]\cap(t_{0}+B_{\varepsilon})}\rho(t)dt.
\end{align}
By (\ref{a}) and  (\ref{b}),  we have
\begin{align*}
\frac{1}{\mu(T,\rho)}\int_{-T}^{T}\left(\int_{t}^{t+1}\mathbb{E}\|h(s)\|^{2}ds\right)^{\frac{1}{2}}\rho(t)dt\geq \frac{\varepsilon \mu(T-\eta,\rho)}{m\mu(T,\rho)\max_{1\leq i\leq m}\{a_{i}\}}
\end{align*}
Since $b=\lim_{T\rightarrow \infty}\sup{\frac{\mu(T-\eta, \rho)}{\mu(T,\rho)}}\in (0, \infty)$, we obtain
\begin{align*}
\lim_{T\rightarrow \infty}\sup\frac{\varepsilon \mu(T-\eta,\rho)}{m\mu(T,\rho)\max_{1\leq i\leq m}\{a_{i}\}}= \frac{b\varepsilon}{m\max_{1\leq i\leq m}\{a_{i}\}}.
\end{align*}
This is contradict with $h^{b}\in SPAA_{0}(\mathbb{R}, L^{2}((0,1), L^{2}(P,H),\rho))$. The proof is complete. \ \ $\Box$\\
\textbf{Theorem 3.4.} Suppose $\rho\in \mathcal{U}_{T}$ and $\lim_{T\rightarrow \infty}\sup{\frac{\mu(T-\eta, \rho)}{\mu(T,\rho)}}>0$ for every $\eta\in \mathbb{R}_{+}$. The space $WPAAS^{2}(\mathbb{R}, L^{2}(P,H),\rho)$ equipped with the norm $\|\cdot\|_{S^{2}}$ is a Banach space. \\
\textbf{Proof.} Assume that $\{f_{n}\}\subset WPAAS^{2}(\mathbb{R},L^{2}(P,H),\rho)$ is a Cauchy sequence with respect to $\|\cdot\|_{S^{2}}$. By definition, there exist $g_{n}$ and $h_{n}$ such that $f_{n}=g_{n}+h_{n}$, where $g_{n}\in AS^{2}(\mathbb{R}, L^{2}(P,H))$ and $h_{n}\in SPAA_{0}(\mathbb{R}, L^{2}((0,1), L^{2}(P,H),\rho))$. By Lemma 3.3, $\{g_{n}(t+\cdot):t\in \mathbb{R}\}\subset \overline{\{f_{n}(t+\cdot):t\in \mathbb{R}\}}$, then $\|g_{n}(t+\cdot)\|_{2}\leq \|f_{n}(t+\cdot)\|_{2}$.  So $\|g_{n}\|_{S^{2}}\leq \|f_{n}\|_{S^{2}}$ and there exist a function $g\in AS^{2}(\mathbb{R}, L^{2}(P,H))$ such that $\|g_{n}-g\|_{S^{2}}\rightarrow 0$ as $n\rightarrow \infty$. Therefore, $\|h_{n}-h\|_{S^{2}}\rightarrow 0$ as $n\rightarrow \infty$, where $h=f-g$. Write $h=(h-h_{n})+h_{n}$, by the Minkowski inequality, we have
\begin{align}\label{limit}
&\frac{1}{\mu(T,\rho)}\int_{-T}^{T}\left(\int_{t}^{t+1}\mathbb{E}\|h(s)\|^{2}\right)^{\frac{1}{2}}\rho(t)dt\nonumber\\
&\leq \frac{1}{\mu(T,\rho)}\int_{-T}^{T}\left(\int_{t}^{t+1}\mathbb{E}\|h_{n}(s)-h(s)\|^{2}\right)^{\frac{1}{2}}\rho(t)dt\nonumber\\
&\quad+\frac{1}{\mu(T,\rho)}\int_{-T}^{T}\left(\int_{t}^{t+1}\mathbb{E}\|h_{n}(s)\|^{2}\right)^{\frac{1}{2}}\rho(t)dt\nonumber\\
&\leq\|h_{n}-h\|_{S^{2}}+\frac{1}{\mu(T,\rho)}\int_{-T}^{T}\left(\int_{t}^{t+1}\mathbb{E}\|h_{n}(s)\|^{2}\right)^{\frac{1}{2}}\rho(t)dt.
\end{align}
Taking $T\rightarrow \infty$ and then taking $n\rightarrow \infty$ with respect to both sides of (\ref{limit}), we have $h^{b}\in SPAA_{0}(\mathbb{R}, L^{2}((0,1), L^{2}(P,H),\rho))$. Since $f_{n}\rightarrow g+h\in WPAAS^{2}(\mathbb{R}, L^{2}(P,H))$ as $n\rightarrow \infty$, we conclude that $WPAAS^{2}(\mathbb{R}, L^{2}(P,H))$ is a Banach space. \ \ $\Box$ \\
\textbf{Definition 3.5.} Let $\rho\in \mathcal{U}_{T}$. A function $f: \mathbb{R}\times L^{2}(P,H)\rightarrow L^{2}(P,H)$, $(t,x)\rightarrow f(t,x)$ with $f(\cdot,x)\in L^{2}_{loc}(\mathbb{R}, L^{2}(P,H))$ for each $x\in L^{2}(P,H)$, is said to be Stepanov-like weighted pseudo almost automorphic (or $S^{2}$-weighted pseudo almost automorphic) if it can be expressed as $f=g+h$, where $g\in AS^{2}(\mathbb{R}\times L^{2}(P,H), L^{2}(P,H))$ and $h^{b}\in SPAA_{0}(\mathbb{R}\times L^{2}(P,H), L^{2}((0,1), L^{2}(P,H)), \rho)$. The collection of such functions will be denoted by  $WPAAS^{2}(\mathbb{R}\times L^{2}(P,H), L^{2}(P,H),\rho)$.\\
\textbf{Lemma 3.6.} Suppose that  $f\in AS^{2}(\mathbb{R}\times L^{2}(P,H))$ and $f(t,u)$ is uniformly stochastically continuous on each bounded subset $K'\subset L^{2}(P,H)$ uniformly for $t\in \mathbb{R}$. If $u\in AS^{2}(L^{2}(P,H))$ and $K=\overline{\{u(t):t\in \mathbb{R}\}}$ is compact. Then $f(\cdot,u(\cdot))\in AS^{2}(L^{2}(P,H))$. \\
\textbf{Proof.} First, we show that $f(\cdot,u(\cdot))\in L^{2}_{loc}(\mathbb{R}, L^{2}(P,H))$. Since $f(t,u)$ is uniformly stochastically continuous on each bounded subset $K'\subset L^{2}(P,H)$ uniformly for $t\in \mathbb{R}$, then for any $\varepsilon>0$, there exists $\delta>0$ such that $x,y\in K'$ and $\mathbb{E}\|x-y\|^{2}< \delta$ imply that
\begin{align}\label{pre}
\mathbb{E}\|f(t,x)-f(t,y)\|^{2}<\varepsilon.
\end{align}
Taking any compact subset $E\subset \mathbb{R}$,
there exists $T_{0}\in \mathbb{N}$ such that $E\subset [-T_{0},T_{0}]$. By (\ref{pre}),
\begin{align}\label{k}
\left(\int_{E}\mathbb{E}\|f(s,x)-f(s,y)\|^{2}ds\right)^{\frac{1}{2}}&\leq \left(\sum_{i=0}^{2T_{0}-1}\int_{-T_{0}+i}^{-T_{0}+i+1}\mathbb{E}\|f(s,x)-f(s,y)\|^{2}ds\right)^{\frac{1}{2}}\nonumber\\
&\leq \sqrt{2T_{0}\varepsilon}
\end{align}
for all $x,y\in K'$ with $\mathbb{E}\|x-y\|^{2}< \delta$.\\
Since $K=\overline{\{u(t):t\in \mathbb{R}\}}$ is compact, there exist finite open balls $O(u_{k},\delta)(k=1,2,\ldots,m)$ defined by
$O(u_{k},\delta)=\{u: \mathbb{E}\|u-u_{k}\|^{2}<\delta\}$ such that $\{u(t):t\in \mathbb{R}\}\subset \bigcup_{k=1}^{m}O(u_{k},\delta)$. Set $B_{k}=\{t\in \mathbb{R}: u(t)\in O(u_{k},\delta)\}$, then $\mathbb{R}=\bigcup _{k=1}^{m}B_{k}$. Set $E_{1}=B_{1}$, $E_{k}=B_{k}\backslash \bigcup _{j=1}^{k-1}B_{j} (2\leq k\leq m)$. Then $E_{i}\cap E_{j}=\emptyset$ for $i\neq j$, $1\leq j\leq m$.\\
Define the step function $\overline{u}: \mathbb{R}\rightarrow L^{2}(P,H)$ by $\overline{u}(s)=u_{k}$, $s\in E_{k}, \ k=1,2,\ldots,m$. It is obvious that $\mathbb{E}\|u(s)-\overline{u}(s)\|^{2}\leq \delta$ for all $s\in \mathbb{R}$. Then we have
\begin{align*}
\left(\int_{E}\mathbb{E}\|f(s,u(s))\|^{2}ds\right)^{\frac{1}{2}}&\leq \left(\int_{E}\mathbb{E}\|f(s,u(s))-f(s,\overline{u}(s))\|^{2}ds\right)^{\frac{1}{2}}+\left(\int_{E}\mathbb{E}\|f(s,\overline{u}(s))\|^{2}ds\right)^{\frac{1}{2}}\\
&\leq \sqrt{2T_{0}\varepsilon}+\left(\sum_{k=1}^{m}\int_{E\cap E_{k}}\mathbb{E}\|f(s,u_{k})\|^{2}\right)^{\frac{1}{2}}.
\end{align*}
Since $f(\cdot,u_{k})\in L^{2}_{loc}(\mathbb{R}, L^{2}(P,H))$, $k=1,2,\ldots, m$. Therefore $f(\cdot,u(\cdot))\in L^{2}_{loc}(\mathbb{R}, L^{2}(P,H))$. \\
Since $f\in AS^{2}(\mathbb{R}\times L^{2}(P,H))$, for every sequence of real numbers $\{s'_{n}\}$, there exist a subsequence $\{s_{n}\}$ and a function $\widetilde{f}: \mathbb{R}\times L^{2}(P,H)\rightarrow L^{2}(P,H)$ with $\widetilde{f}(\cdot,u)\in L^{2}_{loc}(\mathbb{R}, L^{2}(P,H))$ such that
\begin{align}\label{minus}
\int_{t}^{t+1}\mathbb{E}\|f(s+s_{n},u)-\widetilde{f}(s,u)\|^{2}ds\rightarrow 0 \ \mbox{and} \ \int_{t}^{t+1}\mathbb{E}\|\widetilde{f}(s-s_{n},u)-f(s,u)\|^{2}ds\rightarrow 0
\end{align}
as $n\rightarrow \infty$ pointwise on $\mathbb{R}$ and for each $x\in L^{2}(P,H)$. By $u\in AS^{2}(L^{2}(P,H)$, it follows that for every sequence of real numbers $\{s'_{n}\}$, there exist a subsequence $\{s_{n}\}$ and a function $v\in L_{loc}^{2}(\mathbb{R}, L^{2}(P,H))$ such that
\begin{align}\label{compact}
\int_{t}^{t+1}\mathbb{E}\|u(s+s_{n})-v(s)\|^{2}ds\rightarrow 0 \ \mbox{and} \ \int_{t}^{t+1}\mathbb{E}\|v(s-s_{n})-u(s)\|^{2}ds\rightarrow 0
\end{align}
as $n\rightarrow \infty$ pointwise on $\mathbb{R}$.

Since $K=\overline{\{u(t):t\in \mathbb{R}\}}$ is compact, by (\ref{compact}), it follows that $v(t+s)\in K$ for a.e. $s\in [0,1]$. By the assumption of the lemma, for any $\varepsilon>0$, there exists a $\delta>0$ such that $x,y\in K'$ and $\mathbb{E}\|x-y\|^{2}<\delta$ imply that $\mathbb{E}\|f(t,x)-f(t,y)\|^{2}<\varepsilon$ for all $t\in \mathbb{R}$, thus
\begin{align*}
\left(\int_{t}^{t+1}\mathbb{E}\|f(s,x)-f(s,y)\|^{2}ds\right)^{\frac{1}{2}}<\varepsilon
\end{align*}
for each $t\in \mathbb{R}$.\\
By the Minkowski's inequality, we get
\begin{align*}
&\left(\int_{t}^{t+1}\mathbb{E}\|f(s+s_{n},u(s+s_{n})-\widetilde{f}(s,v(s))\|^{2}ds\right)^{\frac{1}{2}}\\
&\leq \left(\int_{t}^{t+1}\mathbb{E}\|f(s+s_{n},u(s+s_{n})-f(s+s_{n},v(s))\|^{2}ds\right)^{\frac{1}{2}}\\
&\quad+\left(\int_{t}^{t+1}\mathbb{E}\|f(s+s_{n},v(s)-\widetilde{f}(s,v(s))\|^{2}ds\right)^{\frac{1}{2}}.
\end{align*}
From the above arguments and (\ref{minus}), it follows that
\begin{align*}
\lim_{n\rightarrow \infty}\left(\int_{t}^{t+1}\mathbb{E}\|f(s+s_{n},u(s+s_{n})-\widetilde{f}(s,v(s))\|^{2}ds\right)^{\frac{1}{2}}=0
\end{align*}
for each $t\in \mathbb{R}$.
Similarly, we can show that
\begin{align*}
\lim_{n\rightarrow \infty}\left(\int_{t}^{t+1}\mathbb{E}\|\widetilde{f}(s-s_{n},v(s-s_{n})-f(s,u(s))\|^{2}ds\right)^{\frac{1}{2}}=0
\end{align*}
for each $t\in \mathbb{R}$.\\
Then $f(\cdot,u(\cdot))\in AS^{2}(L^{2}(P,H))$. The proof is complete. \ \ $\Box$  \\
\textbf{Theorem 3.7. } Let $\rho\in \mathcal{U}_{\infty}$ and let $f=g+h\in WPAAS^{2}(\mathbb{R}\times L^{2}(P,H), L^{2}(P,H),\rho)$ with $g\in AS^{2}(\mathbb{R}\times L^{2}(P,H), L^{2}(P,H))$, $h^{b}\in SPAA_{0}(\mathbb{R}\times L^{2}(P,H), L^{2}((0,1), L^{2}(P,H)), \rho)$. Suppose that
the following conditions $(i)$ and $(ii)$ are satisfied:\\
$(i)$ $f(t,x)$ is Lipschitzian in $x\in L^{2}(P,H)$ uniformly in $t\in \mathbb{R}$, this implies there exists a constant $L>0$ such that
\begin{align*}
\mathbb{E}\|f(t,x)-f(t,y)\|^{2}\leq L\cdot\mathbb{E}\|x-y\|^{2}
\end{align*}
for all $x,y\in L^{2}(P,H)$ and $t\in \mathbb{R}$. \\
$(ii)$ $g(t,x)$ is uniformly continuous in any bounded subset $K'\subset L^{2}(P,H)$ uniformly for $t\in \mathbb{R}$.

If $u=u_{1}+u_{2}\in WPAAS^{2}(L^{2}(P,H),\rho)$, with $u_{1}\in AS^{2}(L^{2}(P,H))$, $u_{2}^{b}\in SPAA_{0}(\mathbb{R}, L^{2}((0,1), L^{2}(P,H),\rho))$ and $K=\overline{\{u_{1}(t):t\in \mathbb{R}\}}$ is compact, then $f(\cdot, u(\cdot))\in WPAAS^{2}(\mathbb{R}, L^{2}(P,H),\rho)$. \\
\textbf{Proof.} By definition, $f$ can be decomposed as
\begin{align*}
f(t,u(t))&=g(t,u_{1}(t))+f(t,u(t))-g(t,u_{1}(t))\\
&=g(t,u_{1}(t))+f(t,u(t))-f(t,u_{1}(t))+h(t,u_{1}(t)).
\end{align*}
Set
\begin{align*}
G(t)=g(t,u_{1}(t)), \ F(t)=f(t,u(t))-f(t,u_{1}(t)), \ H(t)=h(t,u_{1}(t)).
\end{align*}
By Lemma 3.6,  we have $g(\cdot, u_{1}(\cdot))\in AS^{2}(L^{2}(P,H))$.  We will show that $F^{b}(\cdot)+H^{b}(\cdot)\in SPAA_{0}(\mathbb{R}, L^{2}((0,1), L^{2}(P,H),\rho))$. In fact,
\begin{align*}
&\frac{1}{\mu(T,\rho)}\int_{-T}^{T}\left(\int_{t}^{t+1}\mathbb{E}\|f(s,u(s)-f(s,u_{1}(s))\|^{2}ds\right)^{\frac{1}{2}}\rho(t)dt\\
&\leq \frac{L}{\mu(T,\rho)}\int_{-T}^{T}\left(\int_{t}^{t+1}\mathbb{E}\|u(s)-u_{1}(s)\|^{2}ds\right)^{\frac{1}{2}}\rho(t)dt\\
&\leq \frac{L}{\mu(T,\rho)}\int_{-T}^{T}\left(\int_{t}^{t+1}\mathbb{E}\|u_{2}(s)\|^{2}ds\right)^{\frac{1}{2}}\rho(t)dt.
\end{align*}
Since $u_{2}^{b}\in SPAA_{0}(\mathbb{R}, L^{2}((0,1), L^{2}(P,H),\rho))$, then $F^{b}(\cdot)\in SPAA_{0}(\mathbb{R}, L^{2}((0,1), L^{2}(P,H),\rho))$.

Next, we show that $H^{b}(\cdot)\in SPAA_{0}(\mathbb{R}, L^{2}((0,1), L^{2}(P,H),\rho))$. For any $\varepsilon>0$, there exists $\delta>0$ such that
\begin{align}\label{lip}
\mathbb{E}\|g(s,x)-g(s,y)\|^{2}\leq \varepsilon
\end{align}
for all $x,y\in K'$ with $\mathbb{E}\|x-y\|^{2}< \delta$. Set $\delta_{0}=\min\{\delta,\varepsilon\}$, then
\begin{align}\label{sum}
&\left(\int_{t}^{t+1}\mathbb{E}\|h(s,x)-h(s,y)\|^{2}ds\right)^{\frac{1}{2}}\nonumber\\
&\leq \left(\int_{t}^{t+1}\mathbb{E}\|f(s,x)-f(s,y)\|^{2}ds\right)^{\frac{1}{2}}+\left(\int_{t}^{t+1}\mathbb{E}\|g(s,x)-g(s,y)\|^{2}ds\right)^{\frac{1}{2}}\nonumber\\
&\leq \sqrt{L\varepsilon}+\sqrt{\varepsilon}
\end{align}
for $x,y\in K'$ with $\mathbb{E}\|x-y\|^{2}< \delta_{0}$.

By $(ii)$, since $K=\overline{\{u_{1}(t):t\in \mathbb{R}\}}$ is compact, there exist finite open balls $O(x_{k},\delta)(k=1,2,\ldots,m)$ defined by
$O(x_{k},\delta_{0})=\{u: \mathbb{E}\|u-x_{k}\|^{2}<\delta_{0}\}$ such that $\{u_{1}(t):t\in \mathbb{R}\}\subset \bigcup_{k=1}^{m}O(x_{k},\delta_{0})$. \\
Set $B_{k}=\{t\in \mathbb{R}: u_{1}(t)\in O(u_{k},\delta_{0})\}$, then $\mathbb{R}=\bigcup _{k=1}^{m}B_{k}$. Set $E_{1}=B_{1}$, $E_{k}=B_{k}\backslash \bigcup _{j=1}^{k-1}B_{j} (2\leq k\leq m)$. Then $E_{i}\cap E_{j}=\emptyset$ for $i\neq j$, $1\leq j\leq m$. Define the step function $\overline{u}: \mathbb{R}\rightarrow L^{2}(P,H)$ by $\overline{u}(s)=u_{k}$, $s\in E_{k}, \ k=1,2,\ldots,m$. It is obvious that $\mathbb{E}\|u_{1}(s)-\overline{u}(s)\|^{2}\leq \delta_{0}$ for all $s\in \mathbb{R}$. Since $u_{2}^{b}\in SPAA_{0}(\mathbb{R}, L^{2}((0,1), L^{2}(P,H),\rho))$,  there exist $T_{1}>0$ such that for all $T>T_{1}$,
\begin{align}\label{step}
\frac{1}{\mu(T,\rho)}\int_{-T}^{T}\left(\int_{t}^{t+1}\mathbb{E}\|h(s,u_{k})\|^{2}ds\right)^{\frac{1}{2}}\rho(t)dt<\varepsilon, \ k=1,\ldots, m.
\end{align}
By the Minkowski's inequality, for $T>T_{1}$, we have
\begin{align*}
&\frac{1}{\mu(T,\rho)}\int_{-T}^{T}\left(\int_{t}^{t+1}\mathbb{E}\|h(s,u_{1}(s))\|^{2}ds\right)^{\frac{1}{2}}\rho(t)dt\\
&\leq \frac{1}{\mu(T,\rho)}\int_{-T}^{T}\left(\int_{t}^{t+1}\mathbb{E}\|h(s,u_{1}(s))-h(s,\overline{u}(s))\|^{2}ds\right)^{\frac{1}{2}}\rho(t)dt\\
&\quad+\frac{1}{\mu(T,\rho)}\int_{-T}^{T}\left(\int_{t}^{t+1}\mathbb{E}\|h(s,\overline{u}(s))\|^{2}ds\right)^{\frac{1}{2}}\rho(t)dt\\
&\leq\frac{1}{\mu(T,\rho)}\int_{-T}^{T}\left(\int_{t}^{t+1}\mathbb{E}\|h(s,u_{1}(s))-h(s,\overline{u}(s))\|^{2}ds\right)^{\frac{1}{2}}\rho(t)dt\\
&\quad+\frac{1}{\mu(T,\rho)}\int_{-T}^{T}\left(\sum_{k=1}^{m}\int_{[t,t+1]\cap E_{k}}\mathbb{E}\|h(s,u_{k})\|ds\right)^{\frac{1}{2}}\rho(t)dt\\
&\leq\frac{1}{\mu(T,\rho)}\int_{-T}^{T}\left(\int_{t}^{t+1}\mathbb{E}\|h(s,u_{1}(s))-h(s,\overline{u}(s))\|^{2}ds\right)^{\frac{1}{2}}\rho(t)dt\\
&\quad+\frac{1}{\mu(T,\rho)}\sum_{k=1}^{m}\int_{-T}^{T}\left(\int_{[t,t+1]\cap E_{k}}\mathbb{E}\|h(s,u_{k})\|ds\right)^{\frac{1}{2}}\rho(t)dt.
\end{align*}
By (\ref{sum}) and (\ref{step}), we have
\begin{align*}
\frac{1}{\mu(T,\rho)}\int_{-T}^{T}\left(\int_{t}^{t+1}\mathbb{E}\|h(s,u_{1}(s))\|^{2}ds\right)^{\frac{1}{2}}\rho(t)dt\leq \sqrt{L\varepsilon}+\sqrt{\varepsilon}+m\varepsilon.
\end{align*}
Thus $H^{b}(\cdot)\in SPAA_{0}(\mathbb{R}, L^{2}((0,1), L^{2}(P,H),\rho))$, therefore
$f(\cdot, u(\cdot))\in WPAAS^{2}(L^{2}(P,H),\rho)$. The proof is complete. \ \ $\Box$ \\
\section{Main Results}
In this section, we consider the existence and uniqueness of weighted $S^{2}$-pseudo almost automorphic solutions for the stochastic integral equation (\ref{key}).\\
\textbf{Definition 4.1.} (\cite{Pruss}). Let $X$ be a Banach space, $A$ a closed linear unbounded operator in $X$ and $a\in L^{1}_{loc}(\mathbb{R}_{+})$, a scalar kernel $\not\equiv0$.  A family $\{S(t)\}_{t\geq 0}\subset \mathcal{B}(X)$ is called an integral resolvent with generator $A$ if the following conditions are satisfied:

(i) $S(\cdot)x\in L^{1}_{loc}(\mathbb{R}_{+}, X)$ for each $x\in X$ and $\|S(t)\|\leq \phi(t)$ a.e. on $\mathbb{R}_{+}$ for some $\phi\in L^{1}_{loc}(\mathbb{R}_{+})$;

(ii) $S(t)$ commutes with $A$ for each $t\geq 0$;

(iii) the following integral resolvent equation holds
\begin{align*}
S(t)x=a(t)x+\int_{0}^{t}a(t-s)AS(s)xds
\end{align*}
for all $x\in D(A)$ and a.a. $t\geq 0$.

In the following, we recall the concept of $Q$-Wiener process. For more details, we refer to
\cite{LV}. Let $K$ and $H$ be separable Hilbert spaces and $Q$ a self-adjoint nonnegative definite trace-class operator on $K$. The associated eigenvectors forming an orthonormal basis in $K$ will be denoted by $\{e_{k}\}_{k=1}^{\infty}$. Then the space $K_{Q}=Q^{\frac{1}{2}}K$ equipped with the scalar product $\langle u,v\rangle_{K_{Q}}=\sum_{j=1}^{\infty}\frac{1}{\lambda_{j}}\langle u,e_{j}\rangle_{K}\langle v,e_{j}\rangle_{K}$ is a separable Hilbert space with an orthonormal basis $\{\lambda_{j}^{\frac{1}{2}}e_{j}\}_{j=1}^{\infty}$. Denote $L_{2}^{0}=L_{2}^{0}(K_{Q},H)$ the space of Hilbert-Schmidt norm of an operators from $K_{Q}$ to $H$.
Let $\{w_{j}(t)\}_{t\geq 0}$, $j=1,2,\ldots$, be a sequence of independent Brownian motions defined on $(\Omega, \mathcal{F}, P, \mathcal{F}_{t})$. The process $w(t)=\sum_{j=1}^{\infty}\lambda_{j}^{\frac{1}{2}}w_{j}(t)e_{j}$ is called a $Q$-Wiener process in $K$. The two-sided $K$-valued $Q$-Wiener process can be defined as follows: let $\{w_{i}(t), \ t\in \mathbb{R}\}, \ i=1,2$, be independent $K$-valued Wiener processes, then
\begin{equation*}
w(t)=\left\{
\begin{array}{l}
w_{1}(t),\ t\geq 0,\\
w_{2}(-t),\ t\leq 0,
\end{array}
\right.
\end{equation*}
is a two-sided $K$-valued $Q$-Wiener process.  \\
\textbf{Proposition 4.2.} Let $a\in L^{1}(\mathbb{R})$ and $w(t)$ be a two-sided $K$-valued $Q$-wiener process. Suppose that $A$ generates an integral resolvent family $\{S(t)\}_{t\geq 0}$ on $X$,
which satisfies $\|S(t)\|\leq \phi(t)$ for all $t\in \mathbb{R}_{+}$, where $\phi\in L^{1}(\mathbb{R}_{+})$ is nonincreasing. If $f\in L^{1}(\mathbb{R},X)$, $f(t)\in D(A)$ and $\varphi(t)\in D(A)$ P-a.s. for all
$t\in \mathbb{R}$, and
\begin{equation*}
P\left(\int_{-\infty}^{t}\|\varphi(s)\|_{L_{2}^{0}}^{2}ds<\infty\right)=1, \ P\left(\int_{-\infty}^{t}\|A\varphi(s)\|_{L_{2}^{0}}^{2}ds<\infty\right)=1,
\end{equation*}
Then the unique solution of the problem
\begin{align*}
x(t)=\int_{-\infty}^{t}a(t-s)[Ax(s)+f(s)]ds+\int_{-\infty}^{t}a(t-s)\varphi(s)dw(s), \ t\in \mathbb{R},
\end{align*}
is given by
\begin{align*}
x(t)=\int_{-\infty}^{t}S(t-s)f(s)ds+\int_{-\infty}^{t}S(t-s)\varphi(s)dw(s), \ t\in \mathbb{R}.
\end{align*}
\textbf{Proof.} Since $f(t)\in D(A)$ and $A$ is a closed linear operator, from Proposition 1.1.7 in \cite{Arendt}, it follows that $\int_{-\infty}^{t}S(t-s)f(s)ds\in D(A)$. Similar to the proof of Proposition 4.15 in \cite{Da Prato}, we can show that $P\left(\int_{-\infty}^{t}\varphi(s)ds\in D(A)\right)=1$ and $A\int_{-\infty}^{t}\varphi(s)dw(s)=\int_{-\infty}^{t}A\varphi(s)dw(s)$, P-a.s.
By the Fubini theorem, the stochastic Fubini theorem and (iii) of Definition 4.1, we have
\begin{align*}
\int_{-\infty}^{t}a(t-s)Ax(s)ds&=\int_{-\infty}^{t}a(t-s)A\int_{-\infty}^{s}S(s-\tau)f(\tau)d\tau ds\\
&\quad+\int_{-\infty}^{t}a(t-s)A\int_{-\infty}^{s}S(s-\tau)\varphi(\tau)dw(\tau) ds\\
&=\int_{-\infty}^{t}\int_{\tau}^{t}a(t-s)AS(s-\tau)f(\tau)dsd\tau\\
&\quad+\int_{-\infty}^{t}\int_{\tau}^{t}a(t-s)AS(s-\tau)\varphi(\tau)dsdw(\tau)\\
&=\int_{-\infty}^{t}\int_{0}^{t-\tau}a(t-\tau-r)AS(r)drf(\tau)d\tau\\
&\quad+\int_{-\infty}^{t}\int_{0}^{t-\tau}a(t-\tau-r)AS(r)dr\varphi(\tau)dw(\tau)\\
&=\int_{-\infty}^{t}(S(t-\tau)f(\tau)-a(t-\tau)f(\tau))d\tau\\
&\quad+\int_{-\infty}^{t}(S(t-\tau)\varphi(\tau)-a(t-\tau)\varphi(\tau))dw(\tau)\\
&=x(t)-\int_{-\infty}^{t}a(t-\tau)f(\tau)d\tau-\int_{-\infty}^{t}a(t-\tau)\varphi(\tau)dw(\tau).
\end{align*}
The proof is complete. \ \ $\Box$

Motivated by Proposition 4.2, we give the definition of mild solutions for (\ref{key}). \\
\textbf{Definition 4.3.} Let $A$ be the generator of an integral resolvent family $\{S(t)\}_{t\geq 0}$. An $\mathcal{F}_{t}$-progressively measurable stochastic process  $\{x(t)\}_{t\in \mathbb{R}}$ satisfying the stochastic integral equation
\begin{align*}
x(t)=\int_{-\infty}^{t}S(t-s)f(s,x(s))ds+\int_{-\infty}^{t}S(t-s)\varphi(s,x(s))dw(s), \ t\in \mathbb{R},
\end{align*}
is called a mild solution for (\ref{key}).

We list the following assumptions: \\
$(H_{1})$ The operator $A$ is the generator of an integral resolvent family $\{S(t)\}_{t\geq 0}$ on $L^{2}(P,H)$, there exists $\phi\in L^{2}(\mathbb{R}_{+})$, such that $\|S(t)\|\leq \phi(t)$ for all $t\in \mathbb{R}_{+}$. \\
$(H_{2})$ $f, \varphi\in WPAAS^{2}(\mathbb{R}\times L^{2}(P,H), L^{2}(P,H))$, and there exists a constant $L>0$ such that
\begin{align*}
&\mathbb{E}\|f(t,x)-f(t,y)\|^{2}\leq L\cdot\mathbb{E}\|x-y\|^{2}, \\
&\mathbb{E}\|\varphi(t,x)-\varphi(t,y)\|^{2}\leq L\cdot\mathbb{E}\|x-y\|^{2},
\end{align*}
for all $t\in \mathbb{R}$ and each $x,y\in L^{2}(P,H)$. \\
$(H_{3})$ The function $f=g_{1}+h_{1}\in WPAAS^{2}(\mathbb{R}\times L^{2}(P,H), L^{2}(P,H),\rho)$, $\varphi=g_{2}+h_{2}\in WPAAS^{2}(\mathbb{R}\times L^{2}(P,H), L^{2}(P,H),\rho)$, where $g_{1}, g_{2}\in AS^{2}(\mathbb{R}\times L^{2}(P,H), L^{2}(P,H))$ is uniformly stochastically continuous in any bounded subset $M\subset L^{2}(P,H)$ uniformly in $t\in \mathbb{R}$ and $h_{1}^{b}, h_{2}^{b}\in SPAA_{0}(\mathbb{R}\times L^{2}(P,H), L^{2}((0,1), L^{2}(P,H)), \rho)$.\\
\textbf{Lemma 4.4.} Let $S(t)$ be an integral resolvent family satisfying $(H_{1})$, where $\phi: \mathbb{R}_{+}\rightarrow \mathbb{R}_{+}$ is a decreasing function such that $\sum_{n=1}^{\infty}\phi^{2}(n)<\infty$. If $f: \mathbb{R}\rightarrow L^{2}(P,H), \varphi: \mathbb{R}\rightarrow L^{2}(P,H)$ are $S^{2}$-weighted pseudo almost automorphic processes with the weight function $\rho$, and $F(t),\Psi(t)$ are given by
\begin{align*}
F(t)=\int_{-\infty}^{t}S(t-s)f(s)ds, \ \Psi(t)=\int_{-\infty}^{t}S(t-s)\varphi(s)dw(s), \ t\in \mathbb{R},
\end{align*}
then $F,\Psi\in WPAA(\mathbb{R}, L^{2}(P,H),\rho)$. \\
\textbf{Proof.} Since $f,\varphi\in WPAAS^{2}(\mathbb{R},L^{2}(P,H))$, there exist $g_{1},g_{2}\in AS^{2}(\mathbb{R},L^{2}(P,H))$ and $h_{1}^{b},h_{2}^{b}\in SPAA_{0}(\mathbb{R}, L^{2}((0,1), L^{2}(P,H),\rho))$ such that $f=g_{1}+h_{1}, \varphi=g_{2}+h_{2}$.

For $n=1,2,\ldots$, set
\begin{align*}
F_{n}(t)=\int_{n-1}^{n}S(\tau)f(t-\tau)d\tau, \ \Psi_{n}(t)=\int_{n-1}^{n}S(\tau)\varphi(t-\tau)dw(\tau).
\end{align*}

Denote $x_{n}(t)=\int_{n-1}^{n}S(\tau)g_{1}(t-\tau)d\tau$, $y_{n}(t)=\int_{n-1}^{n}S(\tau)g_{2}(t-\tau)dw(\tau)$, $\alpha_{n}(t)=\int_{n-1}^{n}S(\tau)h_{1}(t-\tau)d\tau$, $\beta_{n}(t)=\int_{n-1}^{n}S(\tau)h_{2}(t-\tau)dw(\tau)$, then $F_{n}(t)=x_{n}(t)+\alpha_{n}(t)$, $\Psi_{n}(t)=y_{n}(t)+\beta_{n}(t)$. In order to show that $F_{n}, \Psi_{n}$ are weighted pseudo almost automorphic processes, we need to verify $x_{n},y_{n}\in SAA(\mathbb{R},L^{2}(P,H))$ and $\alpha_{n},\beta_{n}\in SPAA_{0}(\mathbb{R}, L^{2}(P,H),\rho)$. Since
$\left\|\int_{n-1}^{n}S(\tau)g_{1}(s-\tau)d\tau\right\|\leq \int_{n-1}^{n}\|S(\tau)\|\|g_{1}(s-\tau)\|d\tau\leq \phi(n-1)\int_{n-1}^{n}\|g_{1}(s-\tau)\|d\tau$, by the Schwartz inequality, we have
\begin{align*}
\mathbb{E}\|x_{n}(s)\|^{2}&=\mathbb{E}\left\|\int_{n-1}^{n}S(\tau)g_{1}(s-\tau)d\tau\right\|^{2}\\
&\leq \phi^{2}(n-1)\int_{n-1}^{n}\mathbb{E}\left\|g_{1}(s-\tau)\right\|^{2}d\tau\\
&\leq \phi^{2}(n-1)\|g_{1}\|^{2}_{S^{2}}.
\end{align*}
Since $\sum_{n=1}^{\infty}\phi^{2}(n)<\infty$, it follows from the Weirstrass theorem that the sequence of partial sums $\sum_{k=1}^{n}x_{k}(t)$ is uniformly convergent on $\mathbb{R}$. Set $x(t):=\int_{-\infty}^{t}S(t-\tau)g_{1}(\tau)d\tau$. Then $x(t)=\sum_{n=1}^{\infty}x_{n}(t)$. It is obvious that $\|x(t)\|_{2}^{2}\leq \sum_{n=1}^{\infty}\|x_{n}(t)\|_{2}^{2}=\sum_{n=1}^{\infty}\mathbb{E}\|x_{n}(t)\|^{2}\leq \sum_{n=1}^{\infty}\phi^{2}(n-1)\|g_{1}\|^{2}_{S^{2}}$.

Since $g_{1}\in AS^{2}(\mathbb{R},L^{2}(P,H))$, for every sequence of real numbers $\{s_{n}\}$, there exist a subsequence $\{s_{m}\}$ and a function $\widetilde{g_{1}}(\cdot)\in L_{loc}^{2}(\mathbb{R},L^{2}(P,H))$ such that
\begin{align}\label{ad}
\lim_{m\rightarrow\infty}\left(\int_{t}^{t+1}\mathbb{E}\|g_{1}(s+s_{m})-\widetilde{g_{1}}(s)\|^{2}ds\right)^{\frac{1}{2}}=0
\end{align}
and
\begin{align}\label{minu}
\lim_{m\rightarrow\infty}\left(\int_{t}^{t+1}\mathbb{E}\|\widetilde{g_{1}}(s-s_{m})-g_{1}(s)\|^{2}ds\right)^{\frac{1}{2}}=0
\end{align}
hold for each $t\in \mathbb{R}$.

Denote $\widetilde{x}_{n}(t)=\int_{n-1}^{n}S(\tau)g_{1}(t-\tau)d\tau$. By the the Schwartz inequality, we have
\begin{align}\label{doubl}
&\mathbb{E}\|x_{n}(t+s_{m})-\widetilde{x}_{n}(t)\|^{2}\nonumber\\
&=\mathbb{E}\left\|\int_{n-1}^{n}S(\tau)[g_{1}(s+s_{m}-\tau)-\widetilde{g}_{1}(s-\tau)]d\tau\right\|^{2}\nonumber\\
&\leq \phi^{2}(n-1)\int_{n-1}^{n}\mathbb{E}\|g_{1}(s+s_{m}-\tau)-\widetilde{g}_{1}(s-\tau)\|^{2}d\tau.
\end{align}
By (\ref{ad}), (\ref{doubl}), we have
\begin{align*}
\mathbb{E}\|x_{n}(t+s_{m})-\widetilde{x}_{n}(t)\|^{2}\rightarrow 0\ \mbox{as} \ m\rightarrow \infty.
\end{align*}
Similarly, we can show that
\begin{align*}
\mathbb{E}\|\widetilde{x}_{n}(t-s_{m})-x_{n}(t)\|^{2}\rightarrow 0\ \mbox{as} \ m\rightarrow \infty.
\end{align*}
Then $x_{n}(t)\in SAA(\mathbb{R},L^{2}(P,H))$. Hence  $x(t)\in SAA(\mathbb{R},L^{2}(P,H))$.

By the Ito's isometry property of stochastic integral (See \cite{B}, Corollary 3.1.7), we have
\begin{align*}
\mathbb{E}\|y_{n}(t)\|^{2}&=\mathbb{E}\left\|\int_{n-1}^{n}S(\tau)g_{2}(s-\tau)dw(\tau)\right\|^{2}\\
&\leq \int_{n-1}^{n}\|S(\tau)\|^{2}\mathbb{E}\|g_{2}(s-\tau)\|^{2}d\tau\\
&\leq \phi^{2}(n-1)\|g_{2}\|^{2}_{S^{2}}.
\end{align*}
Since $\sum_{n=1}^{\infty}\phi^{2}(n)<\infty$, the sequence of partial sums $\sum_{k=1}^{n}y_{k}(t)$ is uniformly convergent  on $\mathbb{R}$. Set $y(t):=\int_{-\infty}^{t}S(t-\tau)g_{2}(\tau)dw(\tau)$. Then $y(t)=\sum_{n=1}^{\infty}y_{n}(t)$. It is obvious that $\|y(t)\|_{2}^{2}\leq \sum_{n=1}^{\infty}\|y_{n}(t)\|_{2}^{2}=\sum_{n=1}^{\infty}\mathbb{E}\|y_{n}(t)\|^{2}\leq \sum_{n=1}^{\infty}\phi^{2}(n-1)\|g_{2}\|^{2}_{S^{2}}$.

Since $g_{2}\in AS^{2}(\mathbb{R},L^{2}(P,H))$, for every sequence of real numbers $\{s_{n}\}$, there exist a subsequence $\{s_{m}\}$ and a function $\widetilde{g_{2}}(\cdot)\in L_{loc}^{2}(\mathbb{R},L^{2}(P,H))$ such that
\begin{align}\label{add}
\lim_{m\rightarrow\infty}\left(\int_{t}^{t+1}\mathbb{E}\|g_{2}(s+s_{m})-\widetilde{g_{2}}(s)\|^{2}ds\right)^{\frac{1}{2}}=0
\end{align}
and
\begin{align}\label{minus}
\lim_{m\rightarrow\infty}\left(\int_{t}^{t+1}\mathbb{E}\|\widetilde{g_{2}}(s-s_{m})-g_{2}(s)\|^{2}ds\right)^{\frac{1}{2}}=0
\end{align}
hold for each $t\in \mathbb{R}$.

Denote $\widetilde{y}_{n}(t)=\int_{n-1}^{n}S(\tau)g_{2}(t-\tau)dw(\tau)$. By the Ito's isometry property of stochastic integral, we have
\begin{align}\label{double}
&\mathbb{E}\|y_{n}(t+s_{m})-\widetilde{y}_{n}(t)\|^{2}\nonumber\\
&=\mathbb{E}\left\|\int_{n-1}^{n}S(\tau)[g_{2}(s+s_{m}-\tau)-\widetilde{g}_{2}(s-\tau)]dw(\tau)\right\|^{2}\nonumber\\
&\leq \phi^{2}(n-1)\int_{n-1}^{n}\mathbb{E}\|g_{2}(s+s_{m}-\tau)-\widetilde{g}_{2}(s-\tau)\|^{2}d\tau.
\end{align}
By (\ref{add}), (\ref{double}), we have
\begin{align*}
\mathbb{E}\|y_{n}(t+s_{m})-\widetilde{y}_{n}(t)\|^{2}\rightarrow 0\ \mbox{as} \ m\rightarrow \infty.
\end{align*}
Similarly, we can show that
\begin{align*}
\mathbb{E}\|\widetilde{y}_{n}(t-s_{m})-y_{n}(t)\|^{2}\rightarrow 0\ \mbox{as} \ m\rightarrow \infty.
\end{align*}
Then $y_{n}(t)\in SAA(\mathbb{R},L^{2}(P,H))$. Hence  $y(t)\in SAA(\mathbb{R},L^{2}(P,H))$.

Next we prove that $\alpha_{n}(t),\beta_{n}(t)\in SPAA_{0}(\mathbb{R}, L^{2}(P,H),\rho)$. By the Schwartz inequality,
\begin{align*}
\mathbb{E}\|\alpha_{n}(t)\|^{2}&=\mathbb{E}\left\|\int_{n-1}^{n}S(\tau)h_{1}(t-\tau)d\tau\right\|^{2}\\
&\leq \phi^{2}(n-1)\int_{n-1}^{n}\mathbb{E}\left\|h_{1}(t-\tau)\right\|^{2}d\tau\\
\end{align*}
Then we have
\begin{align*}
\frac{1}{\mu(T,\rho)}\int_{-T}^{T}\left(\mathbb{E}\|\alpha_{n}(t)\|^{2}\right)^{\frac{1}{2}}\rho(t)dt&\leq \frac{\phi(n-1)}{\mu(T,\rho)}\int_{-T}^{T}\left(\int_{n-1}^{n}\mathbb{E}\|h_{1}(t-\tau)\|^{2}d\tau\right)^{\frac{1}{2}}\rho(t) dt\\
&\leq \frac{\phi(n-1)}{\mu(T,\rho)}\int_{-T}^{T}\left(\int_{t-n}^{t-n+1}\mathbb{E}\|h_{1}(\tau)\|^{2}d\tau\right)^{\frac{1}{2}}\rho(t)dt.
\end{align*}
Since $h_{1}^{b}\in SPAA_{0}(\mathbb{R}, L^{2}((0,1), L^{2}(P,H),\rho))$, it follows that  $\alpha_{n}(t)\in SPAA_{0}(\mathbb{R}, L^{2}(P,H),\rho)$.
Then the uniform limit $\alpha(t)=\sum_{n=1}^{\infty}\alpha_{n}(t)\in SPAA_{0}(\mathbb{R}, L^{2}(P,H),\rho)$.

By the Ito's isometry property of stochastic integral, we have
\begin{align*}
\mathbb{E}\|\beta_{n}(t)\|^{2}&=\mathbb{E}\left\|\int_{n-1}^{n}S(\tau)h_{2}(t-\tau)dw(\tau)\right\|^{2}\\
&\leq \int_{n-1}^{n}\|S(\tau)\|^{2}\mathbb{E}\|h_{2}(t-\tau)\|^{2}d\tau\\
&\leq \phi^{2}(n-1)\int_{n-1}^{n}\mathbb{E}\|h_{2}(t-\tau)\|^{2}d\tau.
\end{align*}
By the similar argument as above, we can show that $\beta_{n}(t)\in SPAA_{0}(\mathbb{R}, L^{2}(P,H),\rho)$.
Then the uniform limit $\beta(t)=\sum_{n=1}^{\infty}\beta_{n}(t)\in SPAA_{0}(\mathbb{R}, L^{2}(P,H),\rho)$. Hence $F,\Psi\in WPAA(\mathbb{R}, L^{2}(P,H),\rho)$. \ \ $\Box$\\
\textbf{Theorem 4.5.}  Suppose $\rho\in \mathcal{U}^{inv}$ and the hypotheses $(H_{1}), (H_{2}), (H_{3})$ hold. Then the Equation (\ref{key}) has a unique square-mean weighted pseudo almost automorphic mild solution provide that $L\int_{0}^{\infty}\phi^{2}(s)ds<\frac{1}{4}$. \\
\textbf{Proof.} For any $x(t)\in WPAA(\mathbb{R}, L^{2}(P,H))$, define the operator $R$ by
\begin{align*}
(Rx)(t)=\int_{-\infty}^{t}S(t-s)f(s,x(s))ds+\int_{-\infty}^{t}S(t-s)\varphi(s,x(s))dw(s), \ t\in \mathbb{R}.
\end{align*}
Since $\rho\in \mathcal{U}^{inv}$, for any $x(t)\in WPAA(\mathbb{R}, L^{2}(P,H),\rho)$,  by Lemma 3.2, we have  $x(t)\in WPAAS^{2}(\mathbb{R}, L^{2}(P,H),\rho)$. It is easy to show that the range of a stochastically almost automorphic function is relatively compact, then by Theorem 3.7, $f(\cdot,x(\cdot))\in WPAAS^{2}(\mathbb{R},L^{2}(P,H)), \ \varphi(\cdot,x(\cdot))\in WPAAS^{2}(\mathbb{R},L^{2}(P,H))$. By Lemma 4.4, $(Rx)(\cdot)\in WPAA(\mathbb{R},L^{2}(P,H),\rho)$. Therefore $R$ is a mapping from $WPAA(\mathbb{R}, L^{2}(P,H))$ into $WPAA(\mathbb{R}, L^{2}(P,H))$.

For any $x,y\in WPAA(\mathbb{R}, L^{2}(P,H))$,
\begin{align*}
\mathbb{E}\|(Rx)(t)-(Ry)(t)\|^{2}\leq 2\mathbb{E}\|&\int_{-\infty}^{t}S(t-s)[f(s,x(s))-f(s,y(s))]ds\|^{2}\\
&\quad+2\mathbb{E}\int_{-\infty}^{t}S(t-s)[\varphi(s,x(s))-\varphi(s,y(s))]dw(s)\|^{2}\\
&\leq 2\int_{-\infty}^{t}\|S(t-s)\|^{2}\mathbb{E}\|f(s,x(s)-f(s,y(s))\|^{2}ds\\
&\quad+2\mathbb{E}\int_{-\infty}^{t}\|S(t-s)\|^{2}\mathbb{E}\|\varphi(s,x(s))-\varphi(s,y(s))\|^{2}ds\\
&\leq 4L\int_{-\infty}^{t}\|S(t-s)\|^{2}\mathbb{E}\|x(s)-y(s)\|^{2}ds\\
&\leq 4L\int_{0}^{\infty}\|S(s)\|^{2}\mathbb{E}\|x(t-s)-y(t-s)\|^{2}ds\\
&\leq 4L\|x-y\|_{\infty}^{2}\int_{0}^{\infty}\phi^{2}(s)ds.
\end{align*}
Then
\begin{align*}
\|Rx-Ry\|_{\infty}\leq 2\left(L\int_{0}^{\infty}\phi^{2}(s)ds\right)^{\frac{1}{2}}\|x-y\|_{\infty}.
\end{align*}
Then by the Banach fixed point theorem, $R$ has a unique fixed point in $WPAA(\mathbb{R}, L^{2}(P,H))$, which is the unique square-mean weighted pseudo almost automorphic mild solution for the equation (\ref{key}).  The proof is complete. \ \ $\Box$

\end{document}